\documentclass{amsart}

\usepackage{amsmath}
\usepackage{amscd}

\def\HomO{{\mathcal{H}}{\it om}_{{\mathcal{O}}_{X}}}
\def\HomO2{{\mathcal{H}}{\it om}_{{\mathcal{O}}_{X^{2}}}}
\def\HomA{{\mathcal{H}}{\it om}_{{\sf Gr }\mathcal{A}}}
\def\HU{{\underline{{\mathcal{H}}{\it om}}_{{\sf Gr} \mathcal{A}}}}

\def\Hom{{\operatorname{Hom}}}
\def\Homa{{\operatorname{Hom}}_{{\sf Gr}\mathcal{A}}}
\def\Homo{{\operatorname{Hom}}_{{\mathcal{O}}_{X}}}

\def\A{\mathcal{A}}
\def\C{\mathcal{C}}
\def\M{\mathcal{M}}
\def\N{\mathcal{N}}
\def\O{\mathcal{O}}
\def\F{\mathcal{F}}

\def\R{\mathcal{R}}
\def\Q{\mathcal{Q}}

\def\R{\mathcal{R}}
\def\E{\mathcal{E}}

\newtheorem{lemma}{Lemma}[section]

\newtheorem{theorem}[lemma]{Theorem}
\newtheorem{corollary}[lemma]{Corollary}

\theoremstyle{definition}

\newtheorem{definition}[lemma]{\sl Definition}

\theoremstyle{remark}

\begin{document}

\pagenumbering{arabic}

\title{Serre Finiteness and Serre Vanishing For Non-commutative ${\mathbb{P}}^{1}$-bundles}
\author{Adam Nyman}

\thanks{2000 {\it Mathematical Subject Classification. } Primary 14A22; Secondary 16S99}

\begin{abstract}
Suppose $X$ is a smooth projective scheme of finite type over a field $K$, $\mathcal{E}$ is a locally free ${\mathcal{O}}_{X}$-bimodule of rank $2$,  $\mathcal{A}$ is the non-commutative symmetric algebra generated by $\mathcal{E}$ and ${\sf Proj }\A$ is the corresponding non-commutative $\mathbb{P}^{1}$-bundle.  We use the properties of the internal $\operatorname{Hom}$ functor $\HU(-,-)$ to prove versions of Serre finiteness and Serre vanishing for ${\sf Proj }\A$.  As a corollary to Serre finiteness, we prove that ${\sf Proj }\A$ is Ext-finite.  This fact is used in \cite{izu} to prove that if $X$ is a smooth curve over $\operatorname{Spec }K$, ${\sf Proj }\A$ has a Riemann-Roch theorem and an adjunction formula.   
\end{abstract}

\maketitle

{\it Keywords}:  non-commutative geometry, Serre finiteness, non-commutative projective bundle.

\section{introduction}
Non-commutative $\mathbb{P}^{1}$-bundles over curves play a prominent role in the theory of non-commutative surfaces.  For example, certain non-commutative quadrics are isomorphic to non-commutative $\mathbb{P}^{1}$-bundles over curves \cite{quad}.  In addition, every non-commutative deformation of a Hirzebruch surface is given by a non-commutative $\mathbb{P}^{1}$-bundle over $\mathbb{P}^{1}$ \cite[Theorem 7.4.1, p. 29]{p1bundles}.   

The purpose of this paper is to prove versions of Serre finiteness and Serre vanishing (Theorem \ref{theorem.finite} (1) and (2), respectively) for non-commutative ${\mathbb{P}}^{1}$-bundles over smooth projective schemes of finite type over a field $K$.  As a corollary to the first of these results, we prove that such non-commutative ${\mathbb{P}}^{1}$-bundles are $\operatorname{Ext}$-finite.   This fact is used to prove that non-commutative $\mathbb{P}^{1}$-bundles over smooth curves have a Riemann-Roch theorem and an adjunction formula \cite{izu}.  

We now review some important notions from non-commutative algebraic geometry in order to recall the definition of non-commutative ${\mathbb{P}}^{1}$-bundle.  We conclude the introduction by relating the results of this paper to Mori's intersection theory.

If $X$ is a quasi-compact and quasi-separated scheme, then ${\sf Mod }X$, the category of quasi-coherent sheaves on $X$, is a Grothendieck category.  This leads to the following generalization of the notion of scheme, introduced by Van den Bergh in order to define a notion of blowing-up in the non-commutative setting.

\begin{definition} \cite{blowup}
A {\bf quasi-scheme} is a Grothendieck category ${\sf Mod }X$, which we denote by $X$.  $X$ is called a {\bf noetherian} quasi-scheme if the category ${\sf Mod }X$ is locally noetherian.  $X$ is called a {\bf quasi-scheme over $\mathbf{K}$} if the category ${\sf Mod }X$ is $K$-linear.
\end{definition}
If $R$ is a ring and ${\sf Mod }R$ is the category of right $R$-modules, ${\sf Mod }R$ is a quasi-scheme, called the non-commutative affine scheme associated to $R$.  If $A$ is a graded ring, ${\sf Gr }A$ is the category of graded right $A$-modules, ${\sf Tors }A$ is the full subcategory of ${\sf Gr }A$ consisting of direct limits of right bounded modules, and ${\sf Proj }A$ is the quotient category ${\sf Gr }A/{\sf Tors }A$, then ${\sf Proj }A$ is a quasi-scheme called the non-commutative projective scheme associated to $A$.  If $A$ is an Artin-Schelter regular algebra of dimension 3 with the same hilbert series as a polynomial ring in 3 variables, ${\sf Proj }A$ is called a non-commutative ${\mathbb{P}}^{2}$.  

The notion of non-commutative $\mathbb{P}^{1}$-bundle over a smooth scheme $X$ generalizes that of commutative $\mathbb{P}^{1}$-bundle over $X$.  In order to recall the definition of non-commutative $\mathbb{P}^{1}$-bundle, we review some preliminary notions.  Let $S$ be a scheme of finite type over $\operatorname{Spec }K$ and let $X$ be an $S$-scheme.  For $i=1,2$, let $\operatorname{pr}_{i}:X\times_{S}X \rightarrow X$ denote the standard projections, let $\delta:X \rightarrow X\times_{S}X$ denote the diagonal morphism, and let $\Delta$ denote the image of $\delta$.  

\begin{definition}
A {\bf coherent} $\boldsymbol{\mathcal{O}_{X}}${\bf -bimodule}, $\E$, is a coherent $\mathcal{O}_{X \times_{S} X}$-module such that $\operatorname{pr}_{i|\operatorname{Supp \E}}$ is finite for $i=1,2$.  A coherent $\mathcal{O}_{X}$-bimodule $\E$ is {\bf locally free of rank $\boldsymbol{n}$} if $\operatorname{pr}_{i*}\E$ is locally free of rank $n$ for $i=1,2$.
\end{definition} 
Now assume $X$ is smooth.  If $\E$ is a locally free ${\mathcal{O}}_{X}$-bimodule, then let $\E^{*}$ denote the dual of $\E$ \cite[p. 6]{p1bundles}, and let $\E^{j*}$ denote the dual of $\E^{j-1*}$.  Finally, let $\eta: \mathcal{O}_{\Delta} \rightarrow \E \otimes_{\mathcal{O}_{X}} \E^{*}$ denote the counit from $\mathcal{O}_{\Delta}$ to the bimodule tensor product of $\E$ and $\E^{*}$ \cite[p. 7]{p1bundles}.   

\begin{definition} \cite[Section 4.1]{p1bundles}
Let $\E$ be a locally free $\O_{X}$-bimodule.  The {\bf non-commutative symmetric algebra generated by $\boldsymbol{\E}$}, $\mathcal{A}$, is the sheaf-$\mathbb{Z}$-algebra generated by the $\E^{j*}$ subject to the relations $\eta(\mathcal{O}_{\Delta})$.  
\end{definition}
A more explicit definition of non-commutative symmetric algebra is given in Section 2.  We now recall the definition of non-commutative $\mathbb{P}^{1}$-bundle.

\begin{definition} \label{def.pbund} \cite{p1bundles}
Suppose $X$ is a smooth scheme of finite type over $K$, $\E$ is a locally free $\mathcal{O}_{X}$-bimodule of rank 2 and $\A$ is the non-commutative symmetric algebra generated by $\E$.  Let ${\sf Gr} \A$ denote the category of graded right $\A$-modules, let ${\sf Tors }\A$ denote the full subcategory of ${\sf Gr }\A$ consisting of direct limits of right-bounded modules, and let ${\sf Proj }\A$ denote the quotient of ${\sf Gr }\A$ by ${\sf Tors }\A$.  The category ${\sf Proj }\A$ is a {\bf non-commutative $\mathbf{\mathbb{P}^{1}}$-bundle over $\mathbf{X}$.}
\end{definition}
This notion generalizes that of a commutative $\mathbb{P}^{1}$-bundle over $X$ as follows.  Let $\E$ be an $\mathcal{O}_{X}$-bimodule on which $\mathcal{O}_{X}$ acts centrally.  Then $\E$ can be identified with the direct image ${\operatorname{pr}}_{i*}\E$ for $i=1,2$.  If, furthermore, $\E$ is locally free of rank 2 and $\A$ is the non-commutative symmetric algebra generated by $\E$, Van den Bergh proves \cite[Lemma 4.2.1]{p1bundles} that the category ${\sf Proj }\A$ is equivalent to the category ${\sf Mod } \mathbb{P}_{X}(\operatorname{pr}_{i*}\E)$, where $\mathbb{P}_{X}(-)$ is the usual (commutative) projectivization.  

One of the major problems in non-commutative algebraic geometry is to classify non-commutative surfaces.  Since intersection theory on commutative surfaces facilitates the classification of commutative surfaces, one expects intersection theory to be an important tool in non-commutative algebraic geometry.  Mori shows \cite[Theorem 3.11]{izu} that if $Y$ is a noetherian quasi-scheme over a field $K$ such that
\begin{enumerate}
\item{}
$Y$ is $\operatorname{Ext}$-finite,

\item{}
the cohomological dimension of $Y$ is $2$, and

\item{}
$Y$ satisfies Serre duality
\end{enumerate}
then versions of the Riemann-Roch theorem and the adjunction formula hold for $Y$.  Let $X$ be a smooth curve over $\operatorname{Spec }K$.  In \cite{me}, we prove that a non-commutative $\mathbb{P}^{1}$-bundle over $X$ satisfies (2) and (3) above (see Section 4 for a precise statement of these results).  In this paper we prove that a non-commutative $\mathbb{P}^{1}$-bundle over a projective scheme of finite type satisfies (1) (Corollary \ref{corollary.finite}).  We conclude the paper by stating the versions of the Riemann-Roch theorem and the adjunction formula which hold for non-commutative $\mathbb{P}^{1}$-bundles. 
 
In what follows, $K$ is a field, $X$ is a smooth, projective scheme of finite type over $\operatorname{Spec }K$, ${\sf Mod }X$ denotes the category of quasi-coherent $\mathcal{O}_{X}$-modules, and we abuse notation by calling objects in this category $\mathcal{O}_{X}$-modules.
\newline
\newline
{\it Acknowledgment:}  We thank Izuru Mori for showing us his preprint \cite{izu} and for helping us understand the material in Section 4.  

\section{preliminaries}
Before we prove Serre finiteness and Serre vanishing, we review the definition of non-commutative symmetric algebra and the definition and basic properties of the internal Hom functor $\HU(-,-)$ on ${\sf Gr }\A$.

\begin{definition}
Let $\E$ be a locally free $\O_{X}$-bimodule.  The {\bf non-commutative symmetric algebra generated by $\boldsymbol{\E}$} is the sheaf-$\mathbb{Z}$-algebra $\A=\underset{i,j \in \mathbb{Z}}{\oplus}\A_{ij}$ with components
\begin{itemize}
\item{}
$\A_{ii}= \mathcal{O}_{\Delta}$

\item{}
$\A_{i,i+1}=\E^{i*}$,

\item{}
$\A_{ij}= \A_{i,i+1}\otimes \cdots \otimes \A_{j-1,j}/\R_{ij}$ for $j>i+1$, where $\R_{ij} \subset \A_{i,i+1} \otimes \cdots \otimes \A_{j-1,j}$ is the $\mathcal{O}_{X}$-bimodule
$$
\overset{j-2}{\underset{k=i}{\Sigma}}\A_{i,i+1} \otimes \cdots \otimes \A_{k-1,k} \otimes \Q_{k} \otimes \A_{k+2,k+3} \otimes \cdots \otimes \A_{j-1,j},
$$
and $\Q_{i}$ is the image of the unit map $\mathcal{O}_{\Delta} \rightarrow \A_{i,i+1}\otimes \A_{i+1,i+2}$, and
\item{}
$\A_{ij}= 0$ if $i>j$
\end{itemize}
and with multiplication, $\mu$, defined as follows: for $i<j<k$,
\begin{align*}
\A_{ij} \otimes \A_{jk} & = \frac{\A_{i,i+1} \otimes \cdots \otimes \A_{j-1,j}}{\R_{ij}} \otimes \frac{\A_{j,j+1} \otimes \cdots \otimes \A_{k-1,k}}{\R_{jk}} \\
& \cong \frac{\A_{i,i+1} \otimes \cdots \otimes \A_{k-1,k}}{\R_{ij} \otimes \A_{j,j+1} \otimes \cdots \otimes \A_{k-1,k}+ \A_{i,i+1} \otimes \cdots \otimes \A_{j-1,j} \otimes \R_{jk}}
\end{align*}
by \cite[Corollary 3.18]{me2}.  On the other hand,
$$
\R_{ik} \cong \R_{ij} \otimes \A_{j,j+1} \otimes \cdots \otimes \A_{k-1,k}+\A_{i,i+1} \otimes \cdots \otimes \A_{j-1,j} \otimes \R_{jk}+
$$
$$
\A_{i,i+1} \otimes \cdots \otimes \A_{j-2,j-1} \otimes \Q_{j-1} \otimes \A_{j+1,j+2} \otimes \cdots \otimes \A_{k-1,k}.
$$
Thus there is an epi $\mu_{ijk}:\A_{ij} \otimes \A_{jk} \rightarrow \A_{ik}$.  

If $i=j$, let $\mu_{ijk}:\A_{ii} \otimes \A_{ik} \rightarrow \A_{ik}$ be the scalar multiplication map $_{\mathcal{O}}\mu:\mathcal{O}_{\Delta} \otimes \A_{ik} \rightarrow \A_{ik}$.  Similarly, if $j=k$, let $\mu_{ijk}:\A_{ij} \otimes \A_{jj} \rightarrow \A_{ij}$ be the scalar multiplication map $\mu_{\mathcal{O}}$.  Using the fact that the tensor product of bimodules is associative, one can check that multiplication is associative.
\end{definition}

\begin{definition} \label{def.B}
Let $\sf{Bimod }\A-\A$ denote the category of $\A-\A$-bimodules.  Specifically:
\begin{itemize}

\item{}
an object of $\sf{Bimod }\A-\A$ is a triple 
$$
(\mathcal{C}=\{C_{ij}\}_{i,j \in \mathbb{Z}}, \{\mu_{ijk}\}_{i,j,k \in \mathbb{Z}}, \{\psi_{ijk}\}_{i,j,k \in \mathbb{Z}})
$$ 
where ${\mathcal{C}}_{ij}$ is an ${\mathcal{O}}_{X}$-bimodule and $\mu_{ijk}:\C_{ij} \otimes \A_{jk} \rightarrow \C_{ik}$ and $\psi_{ijk}: \A_{ij} \otimes \C_{jk} \rightarrow \C_{ik}$ are morphisms of $\mathcal{O}_{X^{2}}$-modules making $\C$ an $\A$-$\A$ bimodule.

\item{}
A morphism $\phi:  \mathcal{C} \rightarrow \mathcal{D}$ between objects in ${\sf Bimod }\A-\A$ is a  collection $\phi=\{\phi_{ij}\}_{i,j \in \mathbb{Z}}$ such that $\phi_{ij}:{\mathcal{C}}_{ij} \rightarrow {\mathcal{D}}_{ij}$ is a morphism of ${\mathcal{O}}_{X^{2}}$-modules, and such that $\phi$ respects the $\mathcal{A}-\mathcal{A}$-bimodule structure on $\mathcal{C}$ and $\mathcal{D}$.  
\end{itemize}

Let $\mathbb{B}$ denote the full subcategory of $\sf{Bimod }\A-\A$ whose objects $\mathcal{C}=\{C_{ij}\}_{i,j \in \mathbb{Z}}$ have the property that $\C_{ij}$ is coherent and locally free for all $i,j \in \mathbb{Z}$. 

Let ${\mathbb{G}}{\sf r} \mathcal{A}$ denote the full subcategory of $\mathbb{B}$ consisting of objects $\mathcal{C}$ such that for some $n \in \mathbb{Z}$, ${\mathcal{C}}_{ij}=0$ for $i \neq n$ (we say $\C$ is left-concentrated in degree $n$).  
\end{definition}

\begin{definition} \cite[Definition 3.7]{me} \label{def.hom}
Let $\mathcal{C}$ be an object in $\mathbb{B}$ and let $\mathcal{M}$ be a graded right $\mathcal{A}$-module.  We define $\boldsymbol {\HU(\mathcal{C},\mathcal{M})}$ to be the $\mathbb{Z}$-graded $\mathcal{O}_{X}$-module whose $k$th component is the equalizer of the diagram

\begin{equation} \label{eqn.homdef} 
\begin{CD}
\underset{i}{\Pi}{\mathcal{M}}_{i} \otimes {\mathcal{C}}_{ki}^{*} @>\alpha>> \underset{j}{\Pi}{\mathcal{M}}_{j} \otimes {\mathcal{C}}_{kj}^{*}\\
@V{\beta}VV		@VV{\gamma}V\\
\underset{j}{\Pi}(\underset{i}{\Pi}({\mathcal{M}}_{j} \otimes {\mathcal{A}}_{ij}^{*}) \otimes {\mathcal{C}}_{ki}^{*}) @>>{\delta}> \underset{j}{\Pi}(\underset{i}{\Pi}{\mathcal{M}}_{j} \otimes ({\mathcal{C}}_{ki} \otimes {\mathcal{A}}_{ij})^{*}) 
\end{CD}
\end{equation}
where $\alpha$ is the identity map, $\beta$ is induced by the composition 
$$
{\mathcal{M}}_{i} \overset{\eta}{\rightarrow} {\mathcal{M}}_{i}\otimes {\mathcal{A}}_{ij} \otimes {\mathcal{A}}_{ij}^{*} \overset{\mu}{\rightarrow} {\mathcal{M}}_{j}\otimes {\mathcal{A}}_{ij}^{*},
$$
$\gamma$ is induced by the dual of 
$$
{\mathcal{C}}_{ki} \otimes {\mathcal{A}}_{ij} \overset{\mu}{\rightarrow} {\mathcal{C}}_{kj},
$$ 
and $\delta$ is induced by the composition
$$
(\M_{j} \otimes \A_{ij}^{*}) \otimes \C_{ij}^{*} \rightarrow \M_{j}\otimes ({\mathcal{A}}_{ij}^{*}\otimes {\mathcal{C}}_{ki}^{*}) \rightarrow \M_{j} \otimes ({\mathcal{C}}_{ki} \otimes {\mathcal{A}}_{ij})^{*}
$$
whose left arrow is the associativity isomorphism and whose right arrow is induced by the canonical map \cite[Section 2.1]{me}.  If $\mathcal{C}$ is an object of $\mathbb{G}\sf{r} \A$ left-concentrated in degree $k$, we define $\boldsymbol {\HomA (\mathcal{C},\mathcal{M})}$ to be the equalizer of (\ref{eqn.homdef}).
\end{definition}
Let $\tau:{\sf Gr} \A \rightarrow {\sf Tors }\A$ denote the torsion functor, let $\pi: {\sf Gr \A} \rightarrow {\sf Proj }\A$ denote the quotient functor, and let $\omega: {\sf Proj }\A \rightarrow {\sf Gr }\A$ denote the right adjoint to $\pi$.  For any $k \in \mathbb{Z}$, let $e_{k}\A$ denote the right-$\A$-module $\underset{l \geq k}{\bigoplus}\A_{kl}$.  We define $e_{k}\A_{\geq k+n}$ to be the sum $\underset{i \geq 0}{\bigoplus}e_{k}\A_{k+n+i}$ and we let $\A_{\geq n}= \underset{k}{\bigoplus}e_{k}\A_{\geq k+n}$.  

\begin{theorem} \label{theorem.prelim1}
If $\mathcal{M}$ is an object in ${\sf{Gr }}\mathcal{A}$ and $\mathcal{C}$ is an object in $\mathbb{B}$, $\HU (\mathcal{C},\mathcal{M})$ inherits a graded right $\mathcal{A}$-module structure from the left $\mathcal{A}$-module structure of $\mathcal{C}$, making $\HU (-,-):{\mathbb{B}}^{op} \times {\sf{Gr}}\mathcal{A} \rightarrow {\sf{Gr}}\mathcal{A}$ a bifunctor.

Furthermore
\begin{enumerate}
\item{}
$\tau(-) \cong \underset{n \to \infty}{\lim} \HU (\A/\A_{\geq n}, -)$, 

\item{}
If $\F$ is a coherent, locally free $\mathcal{O}_{X}$-bimodule, 
$$
\HomA(\F \otimes e_{k}\A,-) \cong (-)_{k} \otimes \F^{*}
$$ 
and
\item{}
If $\mathcal{L}$ is an $\mathcal{O}_{X}$-module and $\M$ is an object of ${\sf Gr} \A$, 
$$
\Homo(\mathcal{L}, \HomA (e_{k}\mathcal{A}, \mathcal{M})) \cong \operatorname{Hom}_{{\sf{Gr }}\mathcal{A}}(\mathcal{L} \otimes e_{k}\mathcal{A}, \mathcal{M}).
$$
\end{enumerate}
\end{theorem}

\begin{proof}
The first statement is \cite[Proposition 3.11]{me}, (1) is \cite[Proposition 3.19]{me}, (2) is \cite[Theorem 3.16(4)]{me} and (3) is a consequence of \cite[Proposition 3.10]{me} 
\end{proof}
By Theorem \ref{theorem.prelim1} (2), $\HomA(-,\M)$ is $\F \otimes e_{k}\A$-acyclic when $\F$ is a coherent, locally free $\mathcal{O}_{X}$-bimodule.  Thus, one may use the resolution \cite[Theorem 7.1.2]{p1bundles} to compute the derived functors of $\HomA(\A/\A_{\geq 1},-)$.  By Theorem \ref{theorem.prelim1}(1), we may thus compute the derived functors of $\tau$: 
\begin{theorem} \label{theorem.prelim2}
The cohomological dimension of $\tau$ is $2$.  For $i < 2$ and $\mathcal{L}$ a coherent, locally free $\mathcal{O}_{X}$-module,  
$$
\operatorname{R}^{i}\tau(\mathcal{L} \otimes e_{k}\A)=0
$$ 
and 
$$
(\operatorname{R}^{2}\tau(\mathcal{L} \otimes e_{l}\A))_{l-2-i} \cong 
\begin{cases}
\mathcal{L} \otimes \mathcal{Q}_{l-2}^{*} \otimes \A_{l-2-i,l-2}^{*}& \text{if $i \geq 0$}, \\
0& \text{otherwise}.
\end{cases}
$$
\end{theorem}

\begin{proof}
The first result is \cite[Corollary 4.10]{me}, while the remainder is \cite[Lemma 4.9]{me}.
\end{proof}

\section{Serre finiteness and Serre vanishing}
In this section let $I$ denote a finite subset of $\mathbb{Z} \times \mathbb{Z}$.
The proof of the following lemma is straightforward, so we omit it.
\begin{lemma} \label{lemma.coherent}
If $\M$ is a noetherian object in ${\sf Gr} \A$, $\pi \M$ is a noetherian object in ${\sf Proj} \A$ and $\M$ is locally coherent.
\end{lemma}

\begin{lemma} \label{lemma.rtau}
If $\M$ is a noetherian object in ${\sf Gr} \A$, ${\operatorname{R}}^{i}\tau \M$ is locally coherent for all $i \geq 0$.
\end{lemma}

\begin{proof}
The module $\mathcal{O}_{X}(j) \otimes e_{k}\A$ is noetherian by \cite[Lemma 2.17]{me} and the lemma holds with $\M = \underset{(j,k) \in I}{\bigoplus}\mathcal{O}_{X}(j) \otimes e_{k}\A$ by Theorem \ref{theorem.prelim2}.

To prove the result for arbitrary noetherian $\M$, we use descending induction on $i$.  For $i>2$, ${\operatorname{R}}^{i}\tau \M=0$ by Theorem \ref{theorem.prelim2}, so the result is trivial in this case.  Since $\M$ is noetherian, there is a finite subset $I \subset \mathbb{Z} \times \mathbb{Z}$ and a short exact sequence
$$
0 \rightarrow \mathcal{R} \rightarrow {\textstyle \bigoplus\limits_{(j,k) \in I}}\mathcal{O}_{X}(j) \otimes e_{k}\A \rightarrow \M \rightarrow 0
$$
by \cite[Lemma 2.17]{me}.  This induces an exact sequence of $\A$-modules
$$
\ldots \rightarrow ({\operatorname{R}}^{i}\tau (\underset{(j,k) \in I}{\textstyle \bigoplus}\mathcal{O}_{X}(j) \otimes e_{k}\A))_{l} \rightarrow ({\operatorname{R}}^{i}\tau \M)_{l} \rightarrow ({\operatorname{R}}^{i+1}\tau \mathcal{R})_{l} \rightarrow \ldots
$$
The left module is coherent by the first part of the proof, while the right module is coherent by the induction hypothesis.  Hence the middle module is coherent since $X$ is noetherian.
\end{proof}

\begin{corollary} \label{corollary.omega}
If $\M$ is a noetherian object in ${\sf Gr} \A$, ${\operatorname{R}}^{i}(\omega(-)_{k})(\pi \M)$ is coherent for all $i \geq 0$ and all $k \in \mathbb{Z}$.
\end{corollary}

\begin{proof}
Since $(-)_{k}:{\sf Gr }\A \rightarrow {\sf Mod }X$ is an exact functor, ${\operatorname{R}}^{i}(\omega(-)_{k})(\pi \M) \cong {\operatorname{R}}^{i}\omega(\pi \M)_{k}$.  

Now, to prove $\omega(\pi \M)_{k}$ is coherent, we note that there is an exact sequence in ${\sf Mod }X$
$$
0 \rightarrow \tau\M_{k} \rightarrow \M_{k} \rightarrow \omega(\pi \M)_{k} \rightarrow ({\operatorname{R}}^{1}\tau \M)_{k} \rightarrow 0
$$
by \cite[Theorem 4.11]{me}.  Since $\M_{k}$ and $({\operatorname{R}}^{1}\tau \M)_{k}$ are coherent by Lemma \ref{lemma.coherent} and Lemma \ref{lemma.rtau} respectively, $\omega(\pi \M)_{k}$ is coherent since $X$ is noetherian.

The fact that $\operatorname{R}^{i}\omega(\pi \M)_{k}$ is coherent for $i>0$ follows from Lemma \ref{lemma.rtau} since, in this case,
\begin{equation} \label{equation.need}
(\operatorname{R}^{i}\omega(\pi \M))_{k} \cong (\operatorname{R}^{i+1}\tau \M)_{k}
\end{equation}
by \cite[Theorem 4.11]{me}.
\end{proof}

\begin{lemma} \label{lemma.tor}
For $\N$ noetherian in ${\sf Gr }\A$, $\operatorname{R}^{1}\omega(\pi\N)_{k}=0$ for $k>>0$.
\end{lemma}

\begin{proof}
When $\N = \underset{(l,m) \in I}{\textstyle \bigoplus}(\mathcal{O}_{X}(l)\otimes e_{m}\A)$, the result follows from (\ref{equation.need}) and Theorem \ref{theorem.prelim2}.

More generally, there is a short exact sequence  
$$
0 \rightarrow \mathcal{R} \rightarrow \pi (\underset{(l,m) \in I}{\textstyle \bigoplus}\mathcal{O}_{X}(l) \otimes e_{m}\A) \rightarrow \pi\N \rightarrow 0
$$
which induces an exact sequence
$$
\dots \rightarrow \operatorname{R}^{1}\omega(\pi(\underset{(l,m) \in I}{\textstyle \bigoplus}\mathcal{O}_{X}(l) \otimes e_{m}\A)) \rightarrow \operatorname{R}^{1}\omega(\pi\N) \rightarrow \operatorname{R}^{2}\omega(\mathcal{R})=0.
$$
where the right equality is due to (\ref{equation.need}) and Theorem \ref{theorem.prelim2}.  Since the left module is $0$ in high degree, so is $\operatorname{R}^{1}\omega(\pi \N)$.
\end{proof}

\begin{theorem} \label{theorem.finite}
For any noetherian object $\N$ in ${\sf Gr }\A$, 
\begin{enumerate}
\item{}
$\operatorname{Ext}^{i}_{{\sf Proj }\A}(\underset{(j,k) \in I}{\bigoplus}\pi(\mathcal{O}_{X}(j)\otimes e_{k}\A), \pi\N)$ is finite-dimensional over $K$ for all $i \geq 0$, and 
\item{}
for $i > 0$, $\operatorname{Ext}^{i}_{{\sf Proj }\A}(\underset{(j,k) \in I}{\bigoplus}\pi(\mathcal{O}_{X}(j)\otimes e_{k}\A), \pi\N)=0$ whenever $j << 0$ and $k >> 0$.
\end{enumerate}
\end{theorem}

\begin{proof}
Let $d$ denote the cohomological dimension of $X$.  Since $\operatorname{Ext}^{i}_{{\sf Proj }\A}(-, \pi\N)$ commutes with finite direct sums, it suffices to prove the theorem when $I$ has only one element.  
\begin{align*}
\Hom_{{\sf Proj }\A}(\pi(\mathcal{O}_{X}(j) \otimes e_{k}\A),\pi\N) & \cong \Homa(\mathcal{O}_{X}(j) \otimes e_{k}\A, \omega \pi\N) \\
& \cong \Hom_{{\mathcal{O}}_{X}}(\mathcal{O}_{X}(j), \HomA(e_{k}\A,\omega \pi\N)) \\
& \cong \Hom_{{\mathcal{O}}_{X}}(\mathcal{O}_{X}(j), \omega (\pi\N)_{k}) \\
& \cong \Gamma (\mathcal{O}_{X}(-j) \otimes \omega(-)_{k})(\pi\N)
\end{align*}
where the second isomorphism is from Theorem \ref{theorem.prelim1} (3), while the third isomorphism is from Theorem \ref{theorem.prelim1} (2).  Thus, 
$$
\operatorname{Ext}^{i}_{{\sf Proj} \A}(\pi(\mathcal{O}_{X}(j) \otimes e_{k}\A),\pi\N) \cong \operatorname{R}^{i}(\Gamma \circ (\mathcal{O}_{X}(-j) \otimes \omega(-)_{k}))\pi\N.
$$

If $i=0$, (1) follows from Corollary \ref{corollary.omega} and \cite[III, Theorem 5.2a, p. 228]{hart}.  

If $0<i<d+1$, the Grothendieck spectral sequence gives us an exact sequence
\begin{equation} \label{equation.need2}
\ldots \rightarrow \operatorname{R}^{i}\Gamma(\mathcal{O}_{X}(-j) \otimes \omega(\pi\N)_{k}) \rightarrow \operatorname{R}^{i}(\Gamma \circ \mathcal{O}_{X}(-j)\otimes \omega(-)_{k})\pi\N \rightarrow 
\end{equation}
$$
\operatorname{R}^{i-1}\Gamma\operatorname{R}^{1}(\mathcal{O}_{X}(-j)\otimes \omega(-)_{k})\pi\N \rightarrow \ldots 
$$
Since $\omega (\pi\N)_{k}$ and $\operatorname{R}^{1}(\mathcal{O}_{X}(-j)\otimes \omega(-)_{k})\pi\N \cong \mathcal{O}_{X}(-j) \otimes \operatorname{R}^{1}(\omega(-)_{k})\pi\N$ are coherent by Corollary \ref{corollary.omega}, the first and last terms of (\ref{equation.need2}) are finite-dimensional by \cite[III, Theorem 5.2a, p.228]{hart}.  Thus, the middle term of (\ref{equation.need2}) is finite-dimensional as well, which proves (1) in this case.  To prove (2) in this case, we note that, since $\omega(\pi\N)_{k}$ is coherent, the first module of (\ref{equation.need2}) is $0$ for $j<< 0$ by \cite[III, Theorem 5.2b, p.228]{hart}.  If $i>1$, the last module of (\ref{equation.need2}) is $0$ for $j<<0$ for the same reason.  Finally, if $i=1$, the last module of (\ref{equation.need2}) is $0$ since $\operatorname{R}^{1}\omega(\pi\N)_{k}=0$ for $k>>0$ by Lemma \ref{lemma.tor}.

If $i=d+1$, the Grothendieck spectral sequence gives an isomorphism
$$
\operatorname{R}^{d+1}(\Gamma \circ (\mathcal{O}_{X}(-j) \otimes \omega(-)_{k})\pi\N \cong \operatorname{R}^{d}\Gamma\operatorname{R}^{1}(\mathcal{O}_{X}(-j) \otimes \omega(-)_{k})\pi\N. 
$$
In this case, (1) again follows from Corollary \ref{corollary.omega} and \cite[III, Theorem 5.2a, p.228]{hart}, while (2) follows from Lemma \ref{lemma.tor}.
\end{proof}

\begin{corollary} \label{corollary.finite}
If $\M$ and $\N$ are noetherian objects in ${\sf Gr }\A$, $\operatorname{Ext}^{i}_{{\sf Proj} \A}(\pi\M,\pi\N)$ is finite-dimensional for $i \geq 0$.
\end{corollary}

\begin{proof}
Since $\M$ is noetherian, there is an exact sequence
$$
0 \rightarrow \mathcal{R} \rightarrow \pi (\underset{(j,k) \in I}{\textstyle \bigoplus}\mathcal{O}_{X}(j) \otimes e_{k}\A) \rightarrow \pi\M \rightarrow 0.
$$
Since the central term is noetherian by Lemma \ref{lemma.coherent}, so is the $\mathcal{R}$.  Since $\operatorname{Hom}_{{\sf Proj} \A}(-, \pi\N)$ is left exact, there are exact sequences
\begin{equation} \label{eqn.first}
0 \rightarrow \operatorname{Hom}_{{\sf Proj} \A}(\pi\M, \pi\N) \rightarrow \operatorname{Hom}_{{\sf Proj} \A}(\pi (\underset{(j,k) \in I}{\textstyle \bigoplus}\mathcal{O}_{X}(j) \otimes e_{k}\A), \pi\N) \rightarrow
\end{equation}
and, for $i \geq 1$,
\begin{equation} \label{eqn.second}
\rightarrow \operatorname{Ext}^{i-1}_{{\sf Proj} \A}(\mathcal{R},\pi\N) \rightarrow \operatorname{Ext}^{i}_{{\sf Proj} \A}(\pi\M,\pi\N) \rightarrow \operatorname{Ext}^{i}_{{\sf Proj} \A}(\pi (\underset{(j,k) \in I}{\textstyle \bigoplus}\mathcal{O}_{X}(j) \otimes e_{k}\A),\pi\N) \rightarrow 
\end{equation}
Since $\pi$ commutes with direct sums, the right-hand terms of (\ref{eqn.first}) and (\ref{eqn.second}) are finite-dimensional by Theorem \ref{theorem.finite}(1), while the left hand term of (\ref{eqn.second}) is finite-dimensional by the induction hypothesis.
\end{proof}

\section{Riemann-Roch and Adjunction}
Let $X$ be a smooth projective curve, let $\A$ be the noncommutative symmetric algebra generated by a locally free $\mathcal{O}_{X}$-bimodule $\E$ of rank 2, and let $Y = {\sf Proj} \A$.  In this section, we state the Riemann-Roch theorem and adjunction formula for $Y$.  In order to state these results, we need to define an intersection multiplicity on $Y$.  This definition depends on the fact that $Y$ has well behaved cohomology, so we begin this section by reviewing relevant facts regarding the cohomology of $Y$.    

Let $\mathcal{O}_{Y} = \pi \operatorname{pr}_{2*}e_{0}\A$.  By \cite[Theorem 5.20]{me}, $Y$ satisfies Serre duality, i.e., there exists an object $\omega_{Y}$ in ${\sf Proj} \A$, called the canonical sheaf on $Y$, such that 
\begin{equation} \label{eqn.serreduality}
\operatorname{Ext}^{2-i}_{Y}(\mathcal{O}_{Y},-)' \cong \operatorname{Ext}^{i}_{Y}(-,\omega_{Y}) 
\end{equation}
for all $0 \leq i \leq 2$.  Furthermore, the canonical sheaf $\omega_{Y}$ is noetherian \cite{izu0}.

By \cite[Theorem 4.16]{me}, $Y$ has cohomological dimension two, i.e. 
\begin{equation} \label{eqn.cd}
2=\operatorname{sup}\{i|\operatorname{Ext}^{i}_{Y}(\mathcal{O}_{Y},\M) \neq 0 \mbox{ for some noetherian object $\M$ in ${\sf Proj} \A$}\}.
\end{equation}
We write $D:Y \rightarrow Y$ for an autoequivalence, $-D:Y \rightarrow Y$ for the inverse of $D$, and $\M(D):=D(\M)$ for $\M \in Y$.
\begin{definition} \cite[Definition 2.3]{izu}
A {\bf weak divisor} on $Y$ is an element $\mathcal{O}_{D} \in K_{0}(Y)$ of the form $\mathcal{O}_{D}=[\mathcal{O}_{Y}]-[\mathcal{O}_{Y}(-D)]$ for some autoequivalence $D$ of $Y$.
\end{definition}
We now define an intersection multiplicity on $Y$ following \cite{izu}.  Let $\M$ be a noetherian object in ${\sf Proj }\A$, and let $[\M]$ denote its class in $K_{0}(Y)$.  We define, for $\mathcal{O}_{D}$ a weak divisor on $Y$, a map $\xi(\mathcal{O}_{D},-):K_{0}(Y) \rightarrow \mathbb{Z}$ by
$$
\xi(\mathcal{O}_{D},[\M])=\sum\limits_{i=0}^{\infty}(-1)^{i}(\operatorname{dim}_{K}\operatorname{Ext}_{Y}^{i}(\mathcal{O}_{Y},\M)-\operatorname{dim}_{K}\operatorname{Ext}_{Y}^{i}(\mathcal{O}_{Y}(-D),\M)).
$$
This map is well defined by (\ref{eqn.cd}) and Corollary \ref{corollary.finite}.  We define the intersection multiplicity of $\mathcal{O}_{D}$ and $\M$ by
$$
\mathcal{O}_{D} \cdot \M := -\xi(\mathcal{O}_{D},[\M]).
$$
Finally, we define a map $\chi(-):K_{0}(Y) \rightarrow \mathbb{Z}$ by
$$
\chi([\M]):=\sum\limits_{i=0}^{\infty}(-1)^{i}\operatorname{dim}_{K}\operatorname{Ext}_{Y}^{i}(\mathcal{O}_{Y},\M).
$$
\begin{corollary} \label{corollary.rr}
Let $Y = {\sf Proj} \A$, let $\omega_{Y}$ denote the canonical sheaf on $Y$, and suppose $\mathcal{O}_{D}$ is a weak divisor on $Y$.  Then we have the following formulas:
\begin{enumerate}
\item{} (Riemann-Roch)
$$
\chi(\mathcal{O}_{Y}(D)) = \frac{1}{2}(\mathcal{O}_{D} \cdot \mathcal{O}_{D} - \mathcal{O}_{D} \cdot \omega_{Y}+\mathcal{O}_{D} \cdot \mathcal{O}_{Y})+1+p_{a}
$$
where $p_{a} := \chi([\mathcal{O}_{Y}])-1$ is the arithmetic genus of $Y$.

\item{} (Adjunction)
$$
2g-2=\mathcal{O}_{D} \cdot \mathcal{O}_{D}+\mathcal{O}_{D}\cdot \omega_{Y}-\mathcal{O}_{D} \cdot \mathcal{O}_{Y}
$$
where $g:=1-\chi(\mathcal{O}_{D})$ is the genus of $\mathcal{O}_{D}$.
\end{enumerate}
\end{corollary}

\begin{proof}
The quasi-scheme $Y$ is $\operatorname{Ext}$-finite by Corollary \ref{corollary.finite}, has cohomological dimension 2 by \cite[Theorem 4.16]{me}, and satisfies Serre duality with $\omega_{Y}$ by \cite[Theorem 5.20]{me}.  Thus, $Y$ is classical Cohen-Macaulay, and the result follows \cite[Theorem 3.11]{izu}.
\end{proof}    
In stating the Corollary, we defined the intersection multiplicity only for specific elements of $K_{0}(Y) \times K_{0}(Y)$.  In order to define an intersection multiplicity on the entire set $K_{0}(Y) \times K_{0}(Y)$, one must first prove that $Y$ has finite homological dimension.  In \cite[Section 6]{izusmith}, Mori and Smith study noncommutative $\mathbb{P}^{1}$-bundles $Y={\sf Proj} \A$ such that $\A$ is generated by a bimodule $\mathcal{E}$ with the property that $\mathcal{E} \otimes \mathcal{E}$ contains a nondegenerate invertible bimodule.  In this case, they use the structure of $K_{0}(Y)$ to prove that $Y$ has finite homological dimension.  They then compute various intersections on $Y$ without the use of either the Riemann-Roch theorem or the adjunction formula.  In particular, they prove that distinct fibers on $Y$ do not meet, and that a fiber and a section on $Y$ meet exactly once.

\address{Adam Nyman, Department of Mathematical Sciences, Mathematics Building, University of Montana, Missoula, MT 59812-0864}
\end{document}